\documentclass[a4paper,12pt,twoside,leqno,final]{amsart}
\usepackage{amsmath}
\usepackage{amssymb}

\newtheorem{thm}{Theorem}[section]
\newtheorem{lem}[thm]{Lemma}

\newtheorem{prop}[thm]{Proposition}

\newtheorem*{tha}{Theorem A}
\newtheorem*{thb}{Theorem B}

\newcommand{\C}{{\mathbb C}}
\newcommand{\D}{{\mathbb D}}

\newcommand{\T}{{\mathbb T}}

\renewcommand{\sb}{\subset}

\newcommand{\eps}{\varepsilon}

\newcommand{\f}{\frac}
\newcommand{\ov}{\overline}
\newcommand{\al}{\alpha}

\newcommand{\la}{\lambda}
\newcommand{\ze}{\zeta}
\renewcommand{\th}{\theta}

\newcommand{\ph}{\varphi}

\newcommand{\Om}{\Omega}

\newcommand{\const}{\text{\rm const}}

\numberwithin{equation}{section}

\title[Zeros of analytic functions]
{Zeros of analytic functions,\\ 
with or without multiplicities}
\author{Konstantin M. Dyakonov}
\address{ICREA and Universitat de Barcelona, Departament de Matem\`atica 
Aplicada i An\`alisi, Gran Via 585, E-08007 Barcelona, Spain}
\email{dyakonov@mat.ub.es, konstantin.dyakonov@icrea.es}
\thanks{Supported in part by grant MTM2008-05561-C02-01 from El Ministerio de Ciencia 
e Innovaci\'on (Spain) and grant 2009-SGR-1303 from AGAUR (Generalitat de Catalunya).}

\begin{document}
\begin{abstract}
The classical Mason--Stothers theorem deals with nontrivial polynomial 
solutions to the equation $a+b=c$. It provides a lower bound on the number 
of distinct zeros of the polynomial $abc$ in terms of $\deg{a}$, $\deg{b}$ 
and $\deg{c}$. We extend this to general analytic functions living on a reasonable 
bounded domain $\Om\subset\C$, rather than on the whole of $\C$. The estimates 
obtained are sharp, for any $\Om$, and a generalization of the original 
result on polynomials can be recovered from them by a limiting argument. 
\end{abstract}

\maketitle

\section{Introduction}

Given a polynomial $p$ (in one complex variable), write $\deg{p}$ for the degree of $p$ 
and $\widetilde N(p)=\widetilde N_\C(p)$ for the number of its distinct zeros in $\C$. 
The so-called $abc$ theorem, also known as Mason's or Mason--Stothers' theorem, reads as 
follows. 

\begin{tha} Suppose $a$, $b$ and $c$ are relatively prime polynomials, not all constants, 
satisfying $a+b=c$. Then 
\begin{equation}\label{eqn:mason}
\max\{\deg{a},\,\deg{b},\,\deg{c}\}<\widetilde N(abc).
\end{equation}
\end{tha}

Of course, the two sides of \eqref{eqn:mason} are positive integers, so \eqref{eqn:mason} is 
equivalent to saying that the left-hand side does not exceed $\widetilde N(abc)-1$. 

\par More general variants of this result can be found in Mason's book \cite{M}, while 
the current version is essentially due to Stothers \cite{St}. Various approaches to 
and consequences of Theorem A are discussed in \cite{GT, GH, L, ShSm}, the most impressive 
corollary being probably Fermat's Last Theorem for polynomials. The argument leading from $abc$ 
to Fermat is delightfully simple and elegant, so we take the liberty of reproducing it here. 

\par To prove that there are no nontrivial polynomial solutions to the Fermat equation $P^n+Q^n=R^n$ 
for $n\ge3$, apply the $abc$ theorem with $a=P^n$, $b=Q^n$ and $c=R^n$. Write 
$$d=\max\{\deg{P},\,\deg{Q},\,\deg{R}\}$$ 
and note that the left-hand side of \eqref{eqn:mason} equals $nd$. As to the right-hand side, 
$\widetilde N(abc)$, we now have 
$$\widetilde N(abc)=\widetilde N\left((PQR)^n\right)=\widetilde N(PQR)\le\deg{(PQR)}\le3d,$$ 
whence it follows that $n<3$. 

\par The importance of Theorem A is also due to the fact that it served as a prototype (under 
the classical analogy between polynomials and integers) for the famous {\it $abc$ conjecture} 
in number theory. The conjecture, as formulated by Masser and Oesterl\`e in 1985, states that 
to every $\eps>0$ there is a constant $K(\eps)$ with the following property: whenever $a$, $b$ 
and $c$ are relatively prime positive integers satisfying $a+b=c$, one has 
$$c\le K(\eps)\cdot\{\text{\rm rad}(abc)\}^{1+\eps}.$$ 
Here, we write $\text{\rm rad}(\cdot)$ for the {\it radical} of the integer in question, defined 
as the product of the distinct primes that divide it. See, e.\,g., \cite{GT, L} for a discussion 
of the $abc$ conjecture and its potential applications. So far the conjecture remains wide open. 

\par Going back to the polynomial case, let us point out the following \lq\lq$abc\dots xyz$ theorem" 
(or \lq\lq$n$-theorem"), which generalizes Theorem A to sums with any finite number of terms. 

\begin{thb} Let $p_0,p_1,\dots,p_n$ be linearly independent polynomials. Put $p_{n+1}=p_0+\dots+p_n$ 
and assume that the zero-sets $p_0^{-1}(0),\dots,p_{n+1}^{-1}(0)$ are pairwise disjoint. Then 
\begin{equation}\label{eqn:abcxyz}
\max\{\deg{p_0},\dots,\deg{p_{n+1}}\}\le n\widetilde N(p_0p_1\dots p_{n+1})-\f{n(n+1)}2.
\end{equation}
\end{thb}

\par When $n=1$, this reduces to Theorem A. In fact, the assumption on the zero-sets can 
be relaxed to $\bigcap_{j=0}^{n+1} p_j^{-1}(0)=\emptyset$, in which case the quantity 
$\widetilde N(p_0p_1\dots p_{n+1})$ gets replaced by $\sum_{j=0}^{n+1}\widetilde N(p_j)$. 
The latter variant was given by Gundersen and Hayman in \cite[Sect.\,3]{GH}, along with 
a far-reaching generalization from polynomials to entire functions; in addition, 
it was shown there that \eqref{eqn:abcxyz} is asymptotically sharp as $n\to\infty$. 

\par In connection with Theorem B, we also mention Brownawell and Masser's early work 
\cite{BM}, as well as subsequent extensions to polynomials on $\C^m$ (see \cite{HY, SS}) 
and their conjectural analogs in number theory. It should be noted that those extensions -- 
at least the sronger result in \cite{HY} -- relied heavily on Nevanlinna's value 
distribution theory of meromorphic functions, supplemented with some recent developments 
\cite{Ye}. The Nevanlinna theory (or, more precisely, Cartan's version thereof) 
was also the main tool in \cite{GH} when proving the appropriate 
version of Theorem B for entire functions. The analogy between 
Nevanlinna's value distribution theory and Diophantine approximations 
in number theory was unveiled and explored by Vojta \cite{V}. 

\par In this paper, we are concerned with \lq\lq local" versions of the $abc$ theorem -- and 
more generally, of Theorem B -- for analytic functions, a topic not encountered (to the best 
of our knowledge) in the existing literature. This time, the functions will live on a 
bounded -- and reasonably nice -- simply connected domain $\Om\sb\C$ rather than on the whole of 
$\C$. The role of $\deg{p}$, the degree of a polynomial, will of course be played by $N_\Om(f)$, 
the number of zeros (counted with their multiplicities) that $f$ has in $\Om$. Another quantity 
involved will be $\widetilde N_\Om(f)$, the number of the function's distinct zeros in $\Om$. 
The Nevanlinna value distribution theory, which was crucial to earlier approaches in the 
\lq\lq global" setting, will now be replaced by the Riesz--Nevanlinna factorization theory 
on the disk (transplanted, if necessary, to $\Om$). Specifically, Blaschke products will 
be repeatedly employed. 

\par In what follows, we distinguish two cases. First, we assume that the functions involved 
have finitely many zeros in $\Om$ (which enables us to count the zeros, with or without 
multiplicities, and compare the quantities that arise). Of course, this is automatic for 
functions that are analytic on some larger domain containing $\Om\cup\partial\Om$, and we 
actually begin by imposing this stronger assumption. We then weaken the hypotheses, this 
time taking $\Om$ to be the unit disk $\D$, by allowing that the functions be analytic 
on $\D$ and nicely behaved on $\T:=\partial\D$ but still requiring that each of them 
have at most finitely many zeros in $\D$. The results pertaining to this \lq\lq finitely 
many zeros" situation are stated and discussed in Section 2, and then proved in Section 3 
below. 

\par Secondly, we consider the case of infinitely many zeros (the functions being again 
analytic on $\D$ and suitably smooth up to $\T$). Now it makes no sense to count the 
zeros, but the appropriate substitutes for $N_\D(\cdot)$ and $\widetilde N_\D(\cdot)$ 
are introduced and dealt with. This is done in Section 4. 

\par Our method can be roughly described as a mixture of algebraic and analytic techniques. 
The algebraic part is elementary and mimics the reasoning that leads to the classical 
$abc$ theorem, as presented, e.\,g., in \cite{GT}. The analytic component involves certain 
estimates from \cite{C, DHS, VS} that arose when studying the canonical (Riesz--Nevanlinna) 
factorization in various classes of \lq\lq smooth analytic functions". In connection with 
this last topic, which has a long history, let us also mention the seminal paper \cite{H}, 
the monograph \cite{Shi}, and some further developments in \cite{DAmer, DActa, DAdv}. 

\par A special case of our current results from Section 2 was announced in \cite{DCRM}, in 
quite a sketchy form. 

\par We conclude this introduction by asking if our local $abc$-type theorems might suggest, 
by analogy, any number-theoretic results or conjectures. So far, none have occurred to us. 

\section{Finitely many zeros: results and discussion} 

Throughout the rest of the paper, $\Om$ is a bounded simply connected domain in $\C$ such that 
$\partial\Om$ is a rectifiable Jordan curve. We write $dA$ for area measure and $ds$ for arc length, 
and we endow the sets $\Om$ and $\partial\Om$ with the measures $dA/\pi$ and $ds/(2\pi)$, respectively. 
The normalizing factors are chosen so as to ensure that the former (resp., latter) measure 
assigns unit mass to the disk $\D:=\{z\in\C:|z|<1\}$ (resp., to the circle $\T:=\partial\D$). 
The $L^p$-spaces (and norms) on $\Om$ and $\partial\Om$ are then defined in the usual way, with 
respect to the appropriate measure. 

\par In this section, we mainly restrict ourselves to functions that are analytic 
on $\text{\rm clos}\,\Om:=\Om\cup\partial\Om$, i.\,e., analytic on some open set 
containing $\text{\rm clos}\,\Om$. Clearly, such functions (when non-null) can 
only have finitely many zeros in $\Om$, if any. 

\par Let $f$ be analytic on $\text{\rm clos}\,\Om$, and suppose that $a_1,\dots,a_l$ are 
precisely the distinct zeros of $f$ in $\Om$, of multiplicities $m_1,\dots,m_l$ respectively. 
The quantity $N_\Om(f)$, defined as the total number of zeros for $f$ in $\Om$, equals then 
$m_1+\dots+m_l$; the corresponding number of distinct zeros, $\widetilde N_\Om(f)$, is obviously 
$l$. Next, we fix a conformal map $\ph$ from $\Om$ onto $\D$ and write 
\begin{equation}\label{eqn:blaschke}
B(z):=\prod_{k=1}^l\left(\f{\ph(z)-\ph(a_k)}{1-\ov{\ph(a_k)}\ph(z)}\right)^{m_k},
\qquad z\in\Om,
\end{equation}
for the (finite) {\it Blaschke product} built from $f$. The zeros of $B$ in $\Om$, counted with 
multiplicities, are thus the same as those of $f$. In addition, $B$ is continuous up to $\partial\Om$ 
(because $\ph$ is, by Carath\'eodory's theorem) and satisfies $|B(z)|=1$ for all $z\in\partial\Om$. 

\par Now, given Blaschke products $B_1,\dots,B_s$, we write $\text{\rm{LCM}}(B_1,\dots,B_s)$ for 
their {\it least common multiple}, defined in the natural way: this is the Blaschke product 
whose zero-set is $B^{-1}_1(0)\cup\dots\cup B^{-1}_s(0)=:\mathcal Z$, the multiplicity of a zero 
at $a\in\mathcal Z$ being $\max_{1\le j\le s}m(a,B_j)$, where $m(a,B_j)$ is the multiplicity 
of $a$ as a zero of $B_j$. 
\par Further, for a Blaschke product $B$, we 
let $\text{\rm{rad}}(B)$ denote the {\it radical of $B$}; the latter is 
defined (by analogy with the number-theoretic situation) as the Blaschke product 
with zero-set $B^{-1}(0)$ whose zeros are all simple. In other words, if $B$ is given 
by \eqref{eqn:blaschke}, then $\text{\rm{rad}}(B)$ is obtained by replacing each $m_k$ 
with $1$. Observe that $N_\Om(\text{\rm{rad}}(B))=\widetilde N_\Om(B)$. 

\par Finally, we use the notation $W(f_0,\dots,f_n)$ for the {\it Wronskian} of the (analytic) 
functions $f_0,\dots,f_n$, so that 
\begin{equation}\label{eqn:wronskian}
W(f_0,\dots,f_n):=
\begin{vmatrix}
f_0&f_1&\dots&f_n\\
f'_0&f'_1&\dots&f'_n\\
\dots&\dots&\dots&\dots\\
f_0^{(n)}&f_1^{(n)}&\dots&f_n^{(n)}
\end{vmatrix}.
\end{equation}

\par We are now in a position to state the main results of this section. 

\begin{thm}\label{thm:abcdirichlet} Let $f_j$ ($j=0,1,\dots,n$) be analytic 
on $\text{\rm clos}\,\Om$, and suppose that 
the Wronskian $W:=W(f_0,\dots,f_n)$ vanishes nowhere on $\partial\Om$. Let 
\begin{equation}\label{eqn:sum}
f_{n+1}=f_0+\dots+f_n.
\end{equation}
Further, write 
\begin{equation}\label{eqn:lcmrad}
\mathbf B:=\text{\rm{LCM}}(B_0,\dots,B_{n+1})\quad\text{and}\quad
\mathcal B:=\text{\rm{rad}}(B_0B_1\dots B_{n+1}), 
\end{equation}
where $B_j$ is the (finite) Blaschke product associated with $f_j$. Then 
\begin{equation}\label{eqn:maindir}
N_\Om(\mathbf B)\le\la^2+n\mu^2N_\Om(\mathcal B),
\end{equation}
where 
\begin{equation}\label{eqn:defla}
\la=\la_\Om(W):=\|W'\|_{L^2(\Om)}\|1/W\|_{L^\infty(\partial\Om)}
\end{equation}
and 
\begin{equation}\label{eqn:defmu}
\mu=\mu_\Om(W):=\|W\|_{L^\infty(\partial\Om)}\|1/W\|_{L^\infty(\partial\Om)}.
\end{equation}
\end{thm}

In addition to \eqref{eqn:maindir}, we provide an alternative estimate on 
$N_\Om(\mathbf B)$. The factor in front of $N_\Om(\mathcal B)$ will now be 
reduced from $n\mu^2$ to $n\mu$ (note that $\mu\ge1$), while $\la^2$ will 
be replaced by another, possibly larger, quantity. 

\begin{thm}\label{thm:abchardy} Under the hypotheses of Theorem \ref{thm:abcdirichlet}, 
we have 
\begin{equation}\label{eqn:mainhar}
N_\Om(\mathbf B)\le\kappa+n\mu N_\Om(\mathcal B),
\end{equation}
where 
$$\kappa=\kappa_\Om(W):=\|W'\|_{L^1(\partial\Om)}\|1/W\|_{L^\infty(\partial\Om)}$$ 
and $\mu=\mu_\Om(W)$ is defined as in \eqref{eqn:defmu}. 
\end{thm}

\par It should be noted that if the zero-sets of $f_0,\dots,f_{n+1}$ are pairwise disjoint, 
then $\mathbf B=\prod_{j=0}^{n+1}B_j$ and $\mathcal B=\prod_{j=0}^{n+1}\text{\rm rad}(B_j)$, 
in which case 
$$N_\Om(\mathbf B)=\sum_{j=0}^{n+1}N_\Om(f_j)\qquad\text{\rm and}
\qquad N_\Om(\mathcal B)=\sum_{j=0}^{n+1}\widetilde N_\Om(f_j).$$ 

\par The example below shows that the estimates in Theorems \ref{thm:abcdirichlet} 
and \ref{thm:abchardy} are both sharp, for all $\Om$ and $n$, in the sense that equality 
may occur in \eqref{eqn:maindir} and \eqref{eqn:mainhar}. 

\medskip
\noindent\textbf{Example 1.} We may assume, without loss of generality, that $0\in\Om$. Let 
$\Delta$ denote the diameter of $\Om$, and let $\eps$ be a number with $0<\eps<e^{-\Delta}$. 
Consider the functions 
\begin{equation}\label{eqn:fff}
f_0(z)=1,\quad f_j(z)=\eps\f{z^j}{j!}\quad(j=1,\dots,n)
\end{equation}
and set $f_{n+1}=f_0+\dots+f_n$. This done, observe that $f_{n+1}$ has no zeros in $\Om$. Indeed, 
$$\sum_{j=1}^n|f_j(z)|\le\eps\left(e^{|z|}-1\right)\le\eps\left(e^\Delta-1\right)
\le1-\eps,\qquad z\in\Om,$$ 
and so 
$$|f_{n+1}(z)|\ge1-\sum_{j=1}^n|f_j(z)|\ge\eps,\qquad z\in\Om.$$ 
Letting $\ph:\Om\to\D$ be a conformal map with $\ph(0)=0$ (say, the one with $\ph'(0)>0$), 
we see that the Blaschke products $B_j$ associated with the $f_j$'s are given by 
$$B_0(z)=B_{n+1}(z)=1\qquad\text{\rm and}\qquad B_j(z)=\ph^j(z)\quad(j=1,\dots,n).$$ 
Using the notation from \eqref{eqn:lcmrad}, we have then $\mathbf B(z)=\ph^n(z)$ 
and $\mathcal B(z)=\ph(z)$, whence $N_\Om(\mathbf B)=n$ and $N_\Om(\mathcal B)=1$. 
On the other hand, one easily finds that $W:=W(f_0,\dots,f_n)=\eps^n$ (the Wronskian 
matrix being upper triangular), so that $\la_\Om(W)=\kappa_\Om(W)=0$ and $\mu_\Om(W)=1$. 
Consequently, equality holds in both \eqref{eqn:maindir} and \eqref{eqn:mainhar}. 

\medskip The next example, which is a slight modification of the previous one, tells us that 
equality may also occur -- at least on the disk -- in the \lq\lq less trivial" case where 
$W\ne\const$ (or equivalently, when $\la\ne0$ and $\kappa\ne0$). 

\medskip
\noindent\textbf{Example 2.} Let $\Om=\D$ and define $f_0,\dots,f_{n-1}$ as in \eqref{eqn:fff}, 
with some $\eps\in(0,1/e)$. Then put $f_n(z)=\eps z^m/m!$, where $m$ is a fixed integer with 
$m>n$, and finally write $f_{n+1}=f_0+\dots+f_n$. As before, $f_{n+1}$ is zero-free on $\D$. 
For $j=0,\dots,n-1$, the Blaschke product associated with $f_j$ is $z^j$, while 
the Blaschke products corresponding to $f_n$ and $f_{n+1}$ are $z^m$ and $1$, respectively. 
It follows that $\mathbf B(z)=z^m$ and $\mathcal B(z)=z$, whence $N_\D(\mathbf B)=m$ and 
$N_\D(\mathcal B)=1$. The Wronskian $W=W(f_0,\dots,f_n)$ now equals a constant times $z^{m-n}$; 
a simple calculation then yields $\la^2_\D(W)=\kappa_\D(W)=m-n$ and $\mu_\D(W)=1$. Therefore, 
equality is again attained in both \eqref{eqn:maindir} and \eqref{eqn:mainhar}, this time with 
no zero terms on the right. 

\medskip
\par Now let us recall that the functions $f_j$ in Theorems \ref{thm:abcdirichlet} and 
\ref{thm:abchardy} were supposed to be analytic on $\text{\rm clos}\,\Om$. In fact, this can be 
relaxed to the hypothesis that the functions be merely analytic on $\Om$ and suitably smooth 
up to $\partial\Om$. This time, we should explicitly assume that the $f_j$'s have finitely many 
zeros in $\Om$, so as to ensure $N_\Om(\mathbf B)<\infty$ and $N_\Om(\mathcal B)<\infty$. 
The next proposition contains the appropriate versions of the two theorems, specialized 
(for the sake of simplicity) to the case where $\Om=\D$. When stating it, we write 
$\mathcal D=\mathcal D(\D)$ for the {\it Dirichlet space} of the disk, defined as the set 
of all analytic $g$ on $\D$ with $g'\in L^2(\D)$, and we use the standard notation $H^p$ for 
the Hardy spaces on $\D$; see \cite[Chapter II]{G}. 

\begin{prop}\label{prop:dirhardisk} {\rm(a)} Suppose $f_j$ ($j=0,1,\dots,n$) are 
analytic functions on $\D$ satisfying 
\begin{equation}\label{eqn:fnplus1berg}
f_j^{(n)}\in\mathcal D\cap H^\infty
\end{equation}
and $1/W\in L^\infty(\T)$, where $W=W(f_0,\dots,f_n)$. 
Further, put 
$$f_{n+1}=f_0+\dots+f_n$$ 
and assume that 
\begin{equation}\label{eqn:finite}
N_\D(f_j)<\infty\quad\text{for}\quad0\le j\le n+1. 
\end{equation}
Then 
\begin{equation}\label{eqn:lamudisk}
N_\D(\mathbf B)\le\la^2+n\mu^2N_\D(\mathcal B),
\end{equation}
where $\mathbf B$ and $\mathcal B$ are defined as in Theorem \ref{thm:abcdirichlet}, 
$\la=\la_\D(W)$, and $\mu=\mu_\D(W)$. 

\smallskip
{\rm(b)} Replacing \eqref{eqn:fnplus1berg} by the stronger hypothesis that 
\begin{equation}\label{eqn:fnplus1har}
f_j^{(n+1)}\in H^1\qquad
\end{equation}
for all $j$, while retaining the other assumptions above, one has 
\begin{equation}\label{eqn:kamudisk}
N_\D(\mathbf B)\le\kappa+n\mu N_\D(\mathcal B)
\end{equation}
with $\kappa=\kappa_\D(W)$ and $\mu=\mu_\D(W)$. 
\end{prop}

\par We conclude this section by showing that our \lq\lq local" theorems imply Theorem B, 
as stated above, and hence the original $abc$ theorem for polynomials. We shall deduce the 
required \lq\lq global" result from Theorem \ref{thm:abchardy} by a limiting argument. 
An alternative route via Theorem \ref{thm:abcdirichlet} would be equally successful. 

\medskip
\noindent\textbf{Deduction of Theorem B.} 
Suppose $p_0,\dots,p_n$ are linearly independent polynomials and $p_{n+1}=\sum_{j=0}^np_j$. 
Assume also that the zero-sets $p_j^{-1}(0)$ are pairwise disjoint, so that 
$$p_j^{-1}(0)\cap p_k^{-1}(0)=\emptyset\quad\text{\rm whenever}\quad0\le j<k\le n+1.$$ 
An application of Theorem \ref{thm:abchardy} with $\Om=R\D=\{z:|z|<R\}$ gives 
\begin{equation}\label{eqn:manypol}
N_{R\D}(p_0)+\dots+N_{R\D}(p_{n+1})\le\kappa_{R\D}(W)
+n\mu_{R\D}(W)\widetilde N_{R\D}(p_0p_1\dots p_{n+1}).
\end{equation}
Here, 
$$\kappa_{R\D}(W)=\left(\f1{2\pi}\int_{|z|=R}|W'(z)|\,|dz|\right)\cdot
\left(\min_{|z|=R}|W(z)|\right)^{-1}$$ 
and 
$$\mu_{R\D}(W)=\left(\max_{|z|=R}|W(z)|\right)\cdot\left(\min_{|z|=R}|W(z)|\right)^{-1}$$ 
with $W=W(p_0,\dots,p_n)$. Now if $R$ is sufficiently large, then 
$$N_{R\D}(p_j)=N_\C(p_j)=\deg{p_j}=:d_j,\qquad 0\le j\le n+1,$$ 
and 
$$\widetilde N_{R\D}(p_0p_1\dots p_{n+1})=\widetilde N_\C(p_0p_1\dots p_{n+1})=:\widetilde d.$$ 

\par On the other hand, $W$ is a (non-null) polynomial, so that 
$$W(z)=c_mz^m+\text{\rm lower order terms},$$ 
where $m=\deg{W}$ and $c_m\ne0$. The asymptotic behavior of $\kappa_{R\D}(W)$ and $\mu_{R\D}(W)$ 
as $R\to\infty$ is governed by the leading term, $c_mz^m$, whence 
$$\lim_{R\to\infty}\kappa_{R\D}(W)=m\qquad\text{\rm and}\qquad\lim_{R\to\infty}\mu_{R\D}(W)=1.$$ 
We therefore deduce from \eqref{eqn:manypol}, upon letting $R\to\infty$, that 
\begin{equation}\label{eqn:koshka}
d_0+\dots+d_{n+1}\le m+n\widetilde d. 
\end{equation}
To get a bound on $m$, we now recall that $W$ is the sum of $(n+1)!$ products 
of the form 
$$\pm p_0^{(k_0)}p_1^{(k_1)}\dots p_n^{(k_n)},$$ 
where $(k_0,\dots,k_n)$ runs through the permutations of $(0,\dots,n)$. And since 
$$\deg{p_j^{(k_j)}}=d_j-k_j,$$ 
it follows that $m$, the degree of $W$, satisfies 
$$m\le d_0+\dots+d_n-1-\dots-n=d_0+\dots+d_n-\f{n(n+1)}2.$$ 

\par Finally, we put $d:=\max_{0\le j\le n+1}d_j$ and observe that at least two of the 
polynomials involved must be of degree $d$. We may assume that this happens for 
$p_n$ and $p_{n+1}$, so that $d_n=d_{n+1}=d$. The above estimate for $m$ now reads 
$$m\le d_0+\dots+d_{n-1}+d-\f{n(n+1)}2,$$ 
while the left-hand side of \eqref{eqn:koshka} takes the form $d_0+\dots+d_{n-1}+2d$. 
Consequently, \eqref{eqn:koshka} yields 
$$d\le n\widetilde d-\f{n(n+1)}2,$$ 
or equivalently, 
$$\max_{0\le j\le n+1}\deg{p_j}\le n\widetilde N_\C(p_0p_1\dots p_{n+1})-\f{n(n+1)}2,$$ 
as required. 

\section{Finitely many zeros: proofs} 

Let $\mathcal D(\Om)$ denote the {\it Dirichlet space} on $\Om$, i.\,e., the set of all 
analytic functions $f$ on $\Om$ for which the quantity 
\begin{equation}\label{eqn:dirint}
\|f\|^2_{\mathcal D(\Om)}:=\|f'\|^2_{L^2(\Om)}
=\f1\pi\int_\Om|f'(z)|^2dA(z)
\end{equation}
is finite. A bounded analytic function $\th$ on $\Om$ is said to be {\it inner} if its 
nontangential boundary values have modulus $1$ almost everywhere on $\partial\Om$ (with 
respect to arc length). The following result will be needed. 

\begin{lem}\label{lem:factdir} Let $f\in\mathcal D(\Om)$ and let $\th$ be an inner function 
on $\Om$. Then 
\begin{equation}\label{eqn:carldir}
\|f\th\|^2_{\mathcal D(\Om)}=\|f\|^2_{\mathcal D(\Om)}+\f1{2\pi}\int_{\partial\Om}|f|^2|\th'|ds.
\end{equation}
\end{lem} 

In the case where $\Om$ is the unit disk, $\D$, the above lemma follows from Carleson's 
formula in \cite{C}; see also \cite{DHS} for an alternative (operator-theoretic) approach. 
The general case is then established by means of a conformal mapping. Indeed, the Dirichlet 
integral \eqref{eqn:dirint} is conformally invariant, and so is the last term in 
\eqref{eqn:carldir}. 

\par The derivative $\th'$ in \eqref{eqn:carldir} should be interpreted as angular derivative. 
Anyhow, we shall only use formula \eqref{eqn:carldir} when $\th$ is a finite Blaschke product, 
so that $\th=B$ for some $B$ of the form \eqref{eqn:blaschke}. In this situation, $B'$ is sure 
to have nontangential boundary values almost everywhere on $\partial\Om$, since this is the 
case for $\ph'$. Now, applying \eqref{eqn:carldir} to such a $B$ and letting $f\equiv1$, 
we get 
\begin{equation}\label{eqn:bldirhar}
\|B\|^2_{\mathcal D(\Om)}=\f1{2\pi}\int_{\partial\Om}|B'|ds.
\end{equation}
Moreover, the common value of the two sides in \eqref{eqn:bldirhar} is actually $N_\Om(B)$. 
This is clear from the geometric interpretation of the two quantities in terms of area and length, 
combined with the fact that $B$ is an $N$-to-$1$ mapping between $\Om$ and $\D$, where $N=N_\Om(B)$. 

\medskip
\noindent{\it Proof of Theorem \ref{thm:abcdirichlet}.} The first step will be to verify 
that $\mathbf B$ divides $W\mathcal B^n$, in the sense that $W\mathcal B^n/\mathbf B$ is 
analytic on $\Om$. 
\par Clearly, we should only be concerned with those zeros of $\mathbf B$ whose multiplicity 
exceeds $n$. So let $z_0\in\Om$ be a zero of multiplicity $k$, $k>n$, for $\mathbf B$. 
Then there is an index $j\in\{0,\dots,n+1\}$ such that $B_j$ vanishes to order $k$ at $z_0$, 
and so does $f_j$. Expanding the determinant \eqref{eqn:wronskian} along the column that 
contains $f_j,\dots,f_j^{(n)}$, while noting that $f_j^{(l)}$ vanishes to order $k-l$ at $z_0$, 
we see that $W$ has a zero of multiplicity $\ge k-n$ at $z_0$. (In case $j=n+1$, 
one should observe that, by \eqref{eqn:sum}, the determinant remains unchanged upon 
replacing any one of its columns by $(f_{n+1},\dots,f^{(n)}_{n+1})^T$.) And since 
$\mathcal B$ has a zero at $z_0$, it follows that $W\mathcal B^n$ vanishes at least to 
order $k$ at that point. 
\par We conclude that $W\mathcal B^n$ is indeed divisible by $\mathbf B$. In other words, we have 
\begin{equation}\label{eqn:burunduk}
W\mathcal B^n=F\mathbf B, 
\end{equation}
where $F$ is analytic on $\Om$. This $F$ is also continuous on $\text{\rm clos}\,\Om$ because 
$W$, $\mathcal B$ and $\mathbf B$ enjoy this property and because $|\mathbf B|=1$ on $\partial\Om$. 

\par Next, we are going to compute -- and estimate -- the Dirichlet 
integral $\|\cdot\|^2_{\mathcal D(\Om)}$ for each of the two sides of \eqref{eqn:burunduk}. 
On the one hand, an application of Lemma \ref{lem:factdir} yields 
\begin{equation}\label{eqn:estbelow}
\begin{aligned}
\|W\mathcal B^n\|^2_{\mathcal D(\Om)}&=\|F\mathbf B\|^2_{\mathcal D(\Om)}\\
&=\|F\|^2_{\mathcal D(\Om)}+\f1{2\pi}\int_{\partial\Om}|F|^2|\mathbf B'|ds\\
&\ge\f1{2\pi}\int_{\partial\Om}|F|^2|\mathbf B'|ds\\
&\ge\left(\min_{\partial\Om}|F|\right)^2\cdot\f1{2\pi}\int_{\partial\Om}|\mathbf B'|ds\\
&=\|1/W\|_{L^\infty(\partial\Om)}^{-2}N_\Om(\mathbf B).
\end{aligned}
\end{equation}
Here, the last step relies on the fact that $|F|=|W|$ everywhere on $\partial\Om$, an obvious 
consequence of \eqref{eqn:burunduk}. Therefore, the minimum of $|F|$ over $\partial\Om$ coincides 
with that of $|W|$, i.e., with $\|1/W\|_{L^\infty(\partial\Om)}^{-1}$. We have also used the 
equality $(2\pi)^{-1}\int_{\partial\Om}|\mathbf B'|ds=N_\Om(\mathbf B)$, which holds by the 
discussion following \eqref{eqn:bldirhar}. 

\par On the other hand, by Lemma \ref{lem:factdir} again, 
\begin{equation}\label{eqn:estabove}
\begin{aligned}
\|W\mathcal B^n\|^2_{\mathcal D(\Om)}
&=\|W\|^2_{\mathcal D(\Om)}+\f1{2\pi}\int_{\partial\Om}|W|^2|(\mathcal B^n)'|ds\\
&=\|W\|^2_{\mathcal D(\Om)}+\f1{2\pi}\int_{\partial\Om}n|W|^2|\mathcal B'|ds\\
&\le\|W\|^2_{\mathcal D(\Om)}
+n\|W\|_{L^\infty(\partial\Om)}^2\cdot\f1{2\pi}\int_{\partial\Om}|\mathcal B'|ds\\
&=\|W'\|^2_{L^2(\Om)}+n\|W\|_{L^\infty(\partial\Om)}^2N_\Om(\mathcal B).
\end{aligned}
\end{equation} 
Comparing the resulting inequalities from \eqref{eqn:estbelow} and \eqref{eqn:estabove}, we obtain 
$$\|1/W\|_{L^\infty(\partial\Om)}^{-2}N_\Om(\mathbf B)\le\|W'\|^2_{L^2(\Om)}
+n\|W\|_{L^\infty(\partial\Om)}^2N_\Om(\mathcal B),$$ 
which proves \eqref{eqn:maindir}.\quad\qed

To prove Theorem \ref{thm:abchardy}, we need another lemma. Before stating it, we recall 
that an analytic function $f$ on $\Om$ is said to be in the Hardy space $H^p(\Om)$ if 
$(f\circ\psi)\cdot(\psi')^{1/p}$ is in $H^p$ of the disk, for some (or any) conformal map 
$\psi:\D\to\Om$. 

\begin{lem}\label{lem:facthar} Let $f\in H^\infty(\Om)$ and let $\th$ be an inner function 
on $\Om$ with $(f\th)'\in H^1(\Om)$. Then 
\begin{equation}\label{eqn:vinshir}
\|(f\th)'\|_{L^1(\partial\Om)}\ge\f1{2\pi}\int_{\partial\Om}|f|\,|\th'|ds.
\end{equation}
\end{lem} 

For $\Om=\D$, this estimate is due to Vinogradov and Shirokov \cite{VS}. The full statement 
follows by conformal transplantation. Indeed, the class $\{f:f'\in H^1(\Om)\}$ is 
conformally invariant, and so are the two sides of \eqref{eqn:vinshir}. 

\medskip
\noindent{\it Proof of Theorem \ref{thm:abchardy}.} Proceeding as in the proof of Theorem 
\ref{thm:abcdirichlet}, we arrive at \eqref{eqn:burunduk}, where $F$ is analytic on $\Om$ 
and continuous up to $\partial\Om$. Together with Lemma \ref{lem:facthar}, this yields 
\begin{equation}\label{eqn:estbelhar}
\begin{aligned}
\|\left(W\mathcal B^n\right)'\|_{L^1(\partial\Om)}
&=\|\left(F\mathbf B\right)'\|_{L^1(\partial\Om)}\\
&\ge\f1{2\pi}\int_{\partial\Om}|F|\,|\mathbf B'|ds\\
&\ge\left(\min_{\partial\Om}|F|\right)\cdot\f1{2\pi}\int_{\partial\Om}|\mathbf B'|ds\\
&=\|1/W\|_{L^\infty(\partial\Om)}^{-1}N_\Om(\mathbf B).
\end{aligned}
\end{equation}
On the other hand, 
\begin{equation}\label{eqn:estabhar}
\begin{aligned}
\|\left(W\mathcal B^n\right)'\|_{L^1(\partial\Om)}&\le\|W'\mathcal B^n\|_{L^1(\partial\Om)}
+\|W\cdot\left(\mathcal B^n\right)'\|_{L^1(\partial\Om)}\\
&\le\|W'\|_{L^1(\partial\Om)}+n\|W\|_{L^\infty(\partial\Om)}\|\mathcal B'\|_{L^1(\partial\Om)}\\
&=\|W'\|_{L^1(\partial\Om)}+n\|W\|_{L^\infty(\partial\Om)}N_\Om(\mathcal B).
\end{aligned}
\end{equation}
Finally, a juxtaposition of \eqref{eqn:estbelhar} and \eqref{eqn:estabhar} gives 
$$\|1/W\|_{L^\infty(\partial\Om)}^{-1}N_\Om(\mathbf B)\le\|W'\|_{L^1(\partial\Om)}
+n\|W\|_{L^\infty(\partial\Om)}N_\Om(\mathcal B),$$ 
which proves \eqref{eqn:mainhar}.\quad\qed

\medskip
\noindent{\it Proof of Proposition \ref{prop:dirhardisk}.} It is easy to check that 
if either \eqref{eqn:fnplus1berg} or \eqref{eqn:fnplus1har} holds, then the derivatives 
$f_j^{(k)}$ with $0\le k\le n$ are all in $H^\infty$. It follows that, in either case, 
$W\in H^\infty$. Next, note that the derivative $W'$ of the Wronskian $W=W(f_0,\dots,f_n)$ 
is given by 
\begin{equation}\label{eqn:derwron}
W'=
\begin{vmatrix}
f_0&f_1&\dots&f_n\\
f'_0&f'_1&\dots&f'_n\\
\dots&\dots&\dots&\dots\\
f_0^{(n-1)}&f_1^{(n-1)}&\dots&f_n^{(n-1)}\\
f_0^{(n+1)}&f_1^{(n+1)}&\dots&f_n^{(n+1)}
\end{vmatrix}.
\end{equation}
Expanding this determinant along its last row, one therefore deduces that $W\in\mathcal D$ 
in case (a), while $W'\in H^1$ in case (b). These observations show that the quantities 
$\la_\D(W)$, $\mu_\D(W)$ and $\kappa_\D(W)$ appearing in \eqref{eqn:lamudisk} and 
\eqref{eqn:kamudisk} are finite under the stated conditions. 

\par This said, the two estimates are proved in the same way as their counterparts in Theorems 
\ref{thm:abcdirichlet} and \ref{thm:abchardy} above. Namely, one arrives at \eqref{eqn:burunduk} 
as before (with a suitable analytic function $F$ on $\D$) and then essentially rewrites the 
ensuing norm estimates, with $\Om=\D$, based on the (original) disk versions of Lemmas 
\ref{lem:factdir} and \ref{lem:facthar} as contained in \cite{C} and \cite{VS}. 
\par One minor modification is that, in case (a), the functions $W$ and $F$ no longer 
need to be continuous on $\T$. However, they both belong to $\mathcal D\cap H^\infty$, 
and the equality $|F|=|W|$ holds {\it almost everywhere} on $\T$, rather than everywhere. 
Accordingly, the quantity $\min_{\partial\Om}|F|$ appearing in \eqref{eqn:estbelow} 
should be replaced by the essential infimum of $|F|$ over $\T$, which still coincides 
with $\|1/W\|^{-1}_{L^\infty(\T)}$.\quad\qed

\section{Infinitely many zeros}

Given $0<\al\le1$, we write $\mathcal D_\al$ for the space of all analytic functions $f$ 
on $\D$ with 
$$\|f\|^2_{\mathcal D_\al}:=\sum_{k\ge1}k^\al|\widehat f(k)|^2<\infty,$$ 
where $\widehat f(k):=f^{(k)}(0)/k!$. A calculation shows that 
\begin{equation}\label{eqn:equivdiralpha}
\|f\|^2_{\mathcal D_\al}\asymp\f1\pi\int_\D|f'(z)|^2(1-|z|)^{1-\al}dA(z),
\end{equation}
where the notation $U\asymp V$ means that the ratio $U/V$ lies between two positive constants 
depending only on $\al$. When $\al=1$, \eqref{eqn:equivdiralpha} reduces to an identity, and 
$\mathcal D_1$ is just the Dirichlet space $\mathcal D=\mathcal D(\D)$. 

\par Earlier, when proving Theorem \ref{thm:abcdirichlet} and Proposition \ref{prop:dirhardisk} (a), 
we made use of the fact that the total number of zeros of a (finite) Blaschke product $B$ coincides 
with its Dirichlet integral $\|B\|^2_{\mathcal D}$. In this section, we shall be concerned with 
functions living on $\D$ that are allowed to have infinitely many zeros therein. (Our functions 
will, of course, be analytic on $\D$ and appropriately smooth up to $\T$.) The associated Blaschke 
products are thus, in general, infinite products of the form 
\begin{equation}\label{eqn:infbla}
B(z)=z^m\prod_k\left(\f{\bar a_k}{|a_k|}\f{a_k-z}{1-\bar a_kz}\right)^{m_k},\qquad z\in\D; 
\end{equation}
here $a_k$ are the function's distinct zeros in $\D\setminus\{0\}$ of respective 
multiplicities $m_k$, so that $\sum_km_k(1-|a_k|)<\infty$, and $m\ge0$ is the multiplicity 
of its zero at the origin. 

\par While there are no infinite Blaschke products in $\mathcal D=\mathcal D_1$, the spaces 
$\mathcal D_\al$ with $0<\al<1$ do contain such products. (For instance, any Blaschke 
product \eqref{eqn:infbla} satisfying $\sum_km_k(1-|a_k|)^{1-\al}<\infty$ 
will be in $\mathcal D_\al$; see \cite[Theorem 4.2]{A} for this 
and other membership criteria.) Therefore, when looking for 
a reasonable $abc$-type theorem in the current setting, one might expect 
to arrive at a fairly natural formulation by comparing the $\mathcal D_\al$-norms of 
the Blaschke products $\mathbf B$ and $\mathcal B$ rather than counting their zeros. 

\par Here, it is understood that $\mathbf B$ and $\mathcal B$ are built from the given 
functions $f_j$ exactly as before. There is no problem about that, since the notions of 
the least common multiple (LCM) and the radical are perfectly meaningful for infinite 
Blaschke products as well. In particular, if $B$ is defined by \eqref{eqn:infbla} with 
$m\ge0$ and $m_k\ge1$, then $\text{\rm{rad}}(B)$ stands for the product obtained by 
replacing $m$ with $\min\{m,1\}$ and each of the $m_k$'s with $1$. 

\begin{thm}\label{thm:diralpha} Let $0<\al<1$ and suppose $f_j$ ($j=0,1,\dots,n$) are 
analytic functions on $\D$ with 
\begin{equation}\label{eqn:fndiral}
f_j^{(n)}\in\mathcal D_\al\cap H^\infty.
\end{equation}
Assume also that the Wronskian $W:=W(f_0,\dots,f_n)$ satisfies $1/W\in L^\infty(\T)$. 
Put 
$$f_{n+1}=f_0+\dots+f_n.$$ 
Finally, write 
\begin{equation}\label{eqn:lcmradbis}
\mathbf B:=\text{\rm{LCM}}(B_0,\dots,B_{n+1})\quad\text{and}\quad
\mathcal B:=\text{\rm{rad}}(B_0B_1\dots B_{n+1}), 
\end{equation}
where $B_j$ is the Blaschke product associated with $f_j$. Then there exists 
a constant $c_\al>0$ depending only on $\al$ such that 
\begin{equation}\label{eqn:estbis}
c_\al\|\mathbf B\|^2_{\mathcal D_\al}\le\la_\al^2+n\mu^2\|\mathcal B\|^2_{\mathcal D_\al},
\end{equation}
with 
\begin{equation}\label{eqn:deflabis}
\la_\al=\la_{\al,\D}(W):=\|W\|_{\mathcal D_\al}\|1/W\|_\infty
\end{equation}
and 
\begin{equation}\label{eqn:defmubis}
\mu=\mu_\D(W):=\|W\|_\infty\|1/W\|_\infty,
\end{equation}
where $\|\cdot\|_\infty$ stands for $\|\cdot\|_{L^\infty(\T)}$.
\end{thm} 

The proof hinges on the following result, which can be found in \cite[Section 4]{DHS}. 

\begin{lem}\label{lem:fadiral} Let $0<\al<1$. If $f\in\mathcal D_\al$ and $\th$ is an inner 
function on $\D$, then the quantity 
$$\mathcal R_\al(f,\th):=\|f\th\|^2_{\mathcal D_\al}-\|f\|^2_{\mathcal D_\al}$$ 
is nonnegative and satisfies 
$$\mathcal R_\al(f,\th)\asymp\int_\D|f(z)|^2\f{1-|\th(z)|^2}{(1-|z|^2)^{1+\al}}dA(z).$$ 
In particular, 
$$\|\th\|^2_{\mathcal D_\al}\asymp\int_\D\f{1-|\th(z)|^2}{(1-|z|^2)^{1+\al}}dA(z).$$
\end{lem} 

The inequality $\mathcal R_\al(f,\th)\ge0$, when rewritten in the form 
\begin{equation}\label{eqn:divprop}
\|f\th\|_{\mathcal D_\al}\ge\|f\|_{\mathcal D_\al},
\end{equation}
is actually true under the {\it a priori} assumption that $f\in H^2$ and $\th$ is inner. 
This is a refinement of the well-known fact that division by inner factors preserves 
membership in $\mathcal D_\al$; see \cite{H, Shi} and \cite[Section 2]{DAdv} for a discussion 
of a similar phenomenon in various smoothness classes. 

\par Yet another piece of notation will be needed. Namely, given a nonnegative 
measurable function $h$ on $\T$ with $\log h\in L^1(\T)$, we shall write $\mathcal O_h$ 
for the {\it outer function} with modulus $h$, so that 
$$\mathcal O_h(z):=\exp\left(\f1{2\pi}\int_\T\f{\ze+z}{\ze-z}\log h(\ze)\,|d\ze|\right), 
\qquad z\in\D.$$

\medskip
\noindent{\it Proof of Theorem \ref{thm:diralpha}.} First of all, the assumption 
\eqref{eqn:fndiral} implies that $W\in\mathcal D_\al\cap H^\infty$. Indeed, the inclusion 
$W\in H^\infty$ is immediate from the fact that $f_j^{(k)}\in H^\infty$ whenever 
$0\le j,k\le n$. To check that $W\in\mathcal D_\al$, we recall \eqref{eqn:derwron} and 
expand the determinant in that formula along its last row. Since the derivatives $f_j^{(n+1)}$ 
are square integrable against the measure $d\nu_\al(z):=(1-|z|)^{1-\al}dA(z)$, while the 
lower order derivatives are bounded, we infer that $W'\in L^2(d\nu_\al)$ and hence indeed 
$W\in\mathcal D_\al$. Thus, the hypotheses of the theorem guarantee that the quantities 
$\la_\al$ and $\mu$ are finite. We shall also assume that $\mathcal B\in\mathcal D_\al$, 
since otherwise $\|\mathcal B\|_{\mathcal D_\al}=\infty$ and there is nothing to prove. 

\par Arguing as in the proof of Theorem \ref{thm:abcdirichlet}, we verify 
that $\mathbf B$ divides $W\mathcal B^n$, so that \eqref{eqn:burunduk} holds 
with some $F\in H^\infty$. Factoring $F$ canonically (see \cite[Chapter II]{G}), we 
write $F=\mathcal O\mathcal I$, where $\mathcal O$ is outer and $\mathcal I$ is inner. 
Furthermore, since $|F|=|W|$ almost everywhere on $\T$, the outer factor 
$\mathcal O=\mathcal O_{|F|}$ coincides with $\mathcal O_{|W|}$. An application 
of \eqref{eqn:divprop} with $f=\mathcal O_{|W|}\mathbf B$ and $\th=\mathcal I$ 
now shows that 
$$\|W\mathcal B^n\|_{\mathcal D_\al}
=\|F\mathbf B\|_{\mathcal D_\al}
=\|\mathcal O_{|W|}\mathcal I\mathbf B\|_{\mathcal D_\al}
\ge\|\mathcal O_{|W|}\mathbf B\|_{\mathcal D_\al}.$$
From this and Lemma \ref{lem:fadiral} it follows that 
\begin{equation}\label{eqn:wrobelow}
\begin{aligned}
\|W\mathcal B^n\|^2_{\mathcal D_\al}
&\ge\|\mathcal O_{|W|}\mathbf B\|^2_{\mathcal D_\al}\\
&\ge\mathcal R_\al(\mathcal O_{|W|},\mathbf B)\\
&\ge c_1(\al)\int_\D\left|\mathcal O_{|W|}(z)\right|^2
\f{1-|\mathbf B(z)|^2}{(1-|z|^2)^{1+\al}}dA(z)\\
&\ge c_1(\al)\cdot\left(\inf_{z\in\D}\left|\mathcal O_{|W|}(z)\right|\right)^2\cdot
\int_\D\f{1-|\mathbf B(z)|^2}{(1-|z|^2)^{1+\al}}dA(z)\\
&\ge c_2(\al)\cdot\left(\inf_{z\in\D}\left|\mathcal O_{|W|}(z)\right|\right)^2\cdot
\|\mathbf B\|^2_{\mathcal D_\al},
\end{aligned}
\end{equation}
where $c_1(\al)$ and $c_2(\al)$ are positive constants depending only on $\al$. Now let 
us observe that 
$$1/\mathcal O_{|W|}=\mathcal O_{1/|W|}\in H^\infty$$ 
(because $1/W\in L^\infty(\T)$) and 
$$\sup_{z\in\D}\left|\mathcal O_{|W|}(z)\right|^{-1}
=\left\|1/\mathcal O_{|W|}\right\|_\infty=\|1/W\|_\infty,$$ 
whence 
$$\inf_{z\in\D}\left|\mathcal O_{|W|}(z)\right|=\|1/W\|^{-1}_\infty.$$
Substituting this into \eqref{eqn:wrobelow}, we obtain 
\begin{equation}\label{eqn:wrobel}
\|W\mathcal B^n\|^2_{\mathcal D_\al}\ge c_2(\al)\|1/W\|^{-2}_\infty
\|\mathbf B\|^2_{\mathcal D_\al}.
\end{equation}

\par Another application of Lemma \ref{lem:fadiral} in conjunction with the elementary 
inequality $1-t^n\le n(1-t)$, valid for $0\le t\le1$, yields 

\begin{equation*}
\begin{aligned} 
\mathcal R_\al(W,\mathcal B^n)&\le C_1(\al)\int_\D|W(z)|^2\f{1-|\mathcal B(z)|^{2n}}
{(1-|z|^2)^{1+\al}}dA(z)\\
&\le C_1(\al)\cdot n\|W\|^2_\infty\int_\D\f{1-|\mathcal B(z)|^2}
{(1-|z|^2)^{1+\al}}dA(z)\\
&\le C_2(\al)\cdot n\|W\|^2_\infty
\|\mathcal B\|^2_{\mathcal D_\al}, 
\end{aligned}
\end{equation*} 
with suitable constants $C_1(\al)$ and $C_2(\al)$. Consequently, 
\begin{equation}\label{eqn:wroabove}
\|W\mathcal B^n\|^2_{\mathcal D_\al}\le\|W\|^2_{\mathcal D_\al}
+C_2(\al)\cdot n\|W\|^2_\infty\|\mathcal B\|^2_{\mathcal D_\al}.
\end{equation}

\par Finally, a juxtaposition of \eqref{eqn:wrobel} and \eqref{eqn:wroabove} gives 
\begin{equation*}
\begin{aligned} 
c_2(\al)\|1/W\|^{-2}_\infty\|\mathbf B\|^2_{\mathcal D_\al}&\le\|W\|^2_{\mathcal D_\al}
+C_2(\al)\cdot n\|W\|^2_\infty\|\mathcal B\|^2_{\mathcal D_\al}\\
&\le C_2(\al)\cdot\left(\|W\|^2_{\mathcal D_\al}
+n\|W\|^2_\infty\|\mathcal B\|^2_{\mathcal D_\al}\right)
\end{aligned}
\end{equation*}
(we may assume $C_2(\al)\ge1$). This implies \eqref{eqn:estbis}, with 
$$c_\al=c_2(\al)/C_2(\al),$$ 
and completes the proof.\quad\qed

\end{document}